\documentclass[12pt]{article}

\usepackage{amsmath,amssymb,amsfonts,amsthm}

\usepackage{graphicx}
\usepackage{pgf,tikz}
\usetikzlibrary{arrows}

\usepackage{latexsym}
\usepackage{euscript}
\usepackage{mathrsfs}

\usepackage{array}
\usepackage{multirow}

\usepackage[ruled,vlined,Algorithm]{algorithm}
\usepackage{algorithmic}

\usepackage{soul}
\usepackage{verbatim}
\usepackage{alltt}

\usepackage[english]{babel}

\usepackage[numbers]{natbib}

\usepackage{hyperref}
\hypersetup{
  colorlinks=true,
  linkcolor=red,
  citecolor=cyan,
  urlcolor=blue,
  pdfborder={0 0 0}
}

\newtheorem{theorem}{Theorem}[section]

\newtheorem{proposition}{Proposition}[section]

\newtheorem{conclusion}{Conclusion}[section]

\newtheorem{observation}{Observation}[section]

\definecolor{yqqqqq}{rgb}{0.502,0,0}
\definecolor{zzttqq}{rgb}{0.6,0.2,0}
\definecolor{ttqqqq}{rgb}{0.2,0,0}
\definecolor{ccqqww}{rgb}{0.8,0,0.4}

\definecolor{qqqqff}{rgb}{0,0,1}
\definecolor{xdxdff}{rgb}{0.49,0.49,1}
\definecolor{qqffff}{rgb}{0,1,1}

\definecolor{ffqqqq}{rgb}{1,0,0}
\definecolor{ffffqq}{rgb}{1,1,0}
\definecolor{bfffqq}{rgb}{0.749,1,0}

\definecolor{qqffqq}{rgb}{0,1,0}
\definecolor{ccffcc}{rgb}{0.8,1,0.8}

\definecolor{ffffff}{rgb}{1,1,1}
\definecolor{eqeqeq}{rgb}{0.878,0.878,0.878}
\definecolor{cqcqcq}{rgb}{0.753,0.753,0.753}
\definecolor{uuuuuu}{rgb}{0.267,0.267,0.267}

\numberwithin{equation}{section}

\usepackage[a4paper,margin=2.5cm]{geometry}

\begin{document}

\def\nt{\noindent}

\title{A Minimum Doubly Resolving Set and Strong Resolving Set  for the Crystal Cubic Carbon}

\author{Ali Zafari$^{1,}$\footnote{Corresponding author} \and	Saeid Alikhani$^{2}$}

\date{\today}

\maketitle

\begin{center}
$^{1}$Department of Mathematics, Faculty of Science, Payame Noor University, P.O. Box 19395-4697, Tehran, Iran\\ 
{\tt zafari.math@pnu.ac.ir}
\medskip

$^{2}$Department of Mathematical Sciences, Yazd University, 89195-741, Yazd, Iran\\
{\tt alikhani@yazd.ac.ir}
\end{center}
	
\begin{abstract}
The task of identifying resolving sets has been extensively studied due to its wide relevance in fields such as chemistry, robot navigation, combinatorial optimization, pattern recognition, and image processing. These applications have helped motivate and establish the theoretical foundations of the subject. Notably, problems of this type are generally known to be NP-hard.  
This study introduces an alternative structural representation for the crystal cubic carbon \( CCC(n) \). Building on this model, we determine the minimum sizes of both a doubly resolving set and a strong resolving set for $CCC(n)$. 
\end{abstract}


\medskip
\noindent{\bf AMS Subj.\ Class.:}  05C12; 05C90.


\section{Introduction}
\label{sec:introduction}

Throughout this paper, we restrict our attention to connected simple graphs and adopt the notation and terminology of Godsil \cite{Godsil-1}. In graph theory, a graph is commonly described by a set of vertices together with a set of edges connecting them.

The study of resolving sets has a long-standing history in graph theory due to its wide range of applications. These include problems in chemistry \cite{Cameron-1,Chartrand-1}, robot navigation \cite{Khuller-1}, combinatorial optimization \cite{Seb-1}, as well as pattern recognition and image processing, which originally motivated the development of the theory \cite{Slater-1}. Graph theory therefore serves as a fundamental tool for modeling and analyzing chemical structures and complex networks.

In particular, mathematical chemistry-an interdisciplinary field-employs graph-theoretic concepts such as resolving sets, doubly resolving sets, and strong resolving sets to extract structural information. These tools are used to determine the minimum number of atoms required to uniquely identify all other atoms within a molecule. From the perspective of chemical graph theory, a molecular graph represents atoms as vertices and chemical bonds as edges. When a chemical compound is modeled as a graph, variations in atom selection or permutations of their positions generate families of compounds characterized by a shared substructure.

Identifying the smallest set of atoms that uniquely determines the positions of all remaining atoms is a fundamental task in chemistry. Consequently, chemists seek mathematical frameworks that provide distinct representations for distinct compound structures. In graph-theoretic terms, this problem corresponds to determining the minimum cardinality of a resolving set, as discussed in \cite{Baig-1,Erd-1,Khuller-2,Lindst-1,Lindst-2}.

Let $R=\{r_1,r_2,\ldots,r_m\}$ be an ordered subset of vertices of a graph $G$. For any vertex
$u\in V(G)$, the representation of $u$ with respect to $R$ in $G$ is denoted by $r(u\mid R)$
and is defined as the $m$-tuple
\[
r(u\mid R)=(d(u,r_1), d(u,r_2)\ldots, d(u,r_m)),
\]
where $d(u,r_i)$ denotes the length of a shortest path between $u$ and $r_i$ for
$1\le i\le m$. If all vertices in $V(G)\setminus R$ have distinct representations, then
$R$ is called a \emph{resolving set} of $G$ \cite{Buczkowski-1}. Consequently, an important
problem in graph theory is to determine resolving sets of minimum cardinality. The minimum
size of a resolving set of a graph $G$ is called the \emph{metric dimension} of $G$ and is
denoted by $\beta(G)$.

The metric dimension and its related parameters have been extensively studied due to their
numerous applications in graph theory and other scientific disciplines. For further results
and specialized topics, see \cite{Ali-1,Johnson-1,Liu-1,Siddiqui-1,Vetrik-1}. The concept was
first introduced by Slater \cite{Slater-1} under the term \emph{locating set}, while the notion
of metric dimension was independently introduced by Harary and Melter \cite{Harary-1}.

One important extension of the metric dimension is the notion of a \emph{doubly resolving set},
introduced by C\'{a}ceres et al.~\cite{Caceres-1}. An ordered subset
$Z=\{z_1,z_2,\ldots,z_n\}$ of vertices of a graph $G$ is called a doubly resolving set if every
pair of distinct vertices $u,v\in V(G)$ is doubly resolved by some two vertices in $Z$.
Equivalently, for any distinct vertices $u$ and $v$,
\[
r(u\mid Z)-r(v\mid Z)\neq \mu I,
\]
where $\mu$ is an integer and $I=(1,\ldots,1)$ denotes the unit $n$-vector. The minimum
cardinality of a doubly resolving set of $G$ is denoted by $\psi(G)$. Further information on
doubly resolving set

\section{Structure of $CCC(n)$ }
\label{sec:introduction}

The molecular graph of crystal cubic carbon $CCC(2)$ is illustrated in Figure~3. The construction of this structure begins with a single unit cube, after which additional cubes are attached at each vertex of the initial cube through edges. Further details regarding the structural properties of crystal cubic carbon $CCC(n)$ can be found in \cite{Baig-1}. It is well known that the numbers of vertices and edges of $CCC(n)$ are given, respectively, by
\begin{align*}
	|V(\text{CCC}(n))| &= 2 \left\{ 24 \sum_{r=3}^{n} (2^3 - 1)^{r-3} + 31(2^3 - 1)^{n-2} + 2 \sum_{r=0}^{n-2} (2^3 - 1)^r + 3 \right\},  \\
	|E(\text{CCC}(n))| &= 4 \left\{ 24 \sum_{r=3}^{n} (2^3 - 1)^{r-3} + 24(2^3 - 1)^{n-2} + 2 \sum_{r=0}^{n-2} (2^3 - 1)^r + 3 \right\}.
\end{align*}

Several studies have investigated the topological characteristics of the graph $CCC(n)$. Among the most notable degree-based topological indices computed for this structure are the Randi\'c, ABC, and Zagreb indices, along with other related invariants; see \cite{Gao-1,Yang-1}.

In \cite{Zhang-1}, the authors determined the metric dimension of $CCC(n)$. Specifically, Theorem~1 in \cite{Zhang-1} establishes that for any integer $n\geq 2$, the minimum cardinality of a resolving set of $CCC(n)$ is equal to $16\times 7^{(n-2)}$. The primary objective of the present work is to compute the minimum cardinalities of both doubly resolving sets and strong resolving sets of the crystal cubic carbon graph $CCC(n)$. To achieve this goal, we first introduce an alternative description of the structure of $CCC(n)$, which facilitates a systematic relabeling of all vertices and is essential for determining these resolving parameters.

\medskip

\noindent
(1) Let $n$ be a fixed positive integer and let $k$ be an integer satisfying $2\leq k\leq n$. Define the sets
$W_1=\{1,\ldots,4\}$ and $W_2=\{5,\ldots,8\}$, and let $L_1$ be a graph with vertex set $\{1,2,\ldots,8\}$. The edge set $E(L_1)$ is defined by
\begin{align*}
	E(L_1)
	&=
	\bigcup_{k=1}^{2}
	\{\, ij \mid i,j \in W_k,\; i<j,\; j-i \in \{1,3\} \,\}
	\;\cup\;
	\{\, ij \mid i \in W_1,\; j \in W_2,\; j-i=4 \,\}.
\end{align*}
It can be verified that the graph $L_1$ is isomorphic to the Cartesian product $C_4 \Box P_2$, where $C_4$ and $P_2$ denote the cycle on four vertices and the path on two vertices, respectively. Henceforth, for simplicity, we identify $V(C_4 \Box P_2)$ with $V(L_1)$ and $E(C_4 \Box P_2)$ with $E(L_1)$. The graph $C_4 \Box P_2$ is shown in Figure~1.

\begin{center}
\begin{tikzpicture}[line cap=round,line join=round,>=triangle 45,x=15.0cm,y=10.0cm]
\clip(5.531538312092156,2.3252649721940215) rectangle (6.113499725856829,2.865884864762297);
\draw (5.999118452553847,2.7728532194486912)-- (6.,2.6);
\draw (6.,2.6)-- (5.8,2.6);
\draw (5.8,2.6)-- (5.797172680635512,2.766733650602681);
\draw (5.797172680635512,2.766733650602681)-- (5.999118452553847,2.7728532194486912);
\draw (5.859090364648229,2.6814380778860616)-- (5.865450598678444,2.512891876085364);
\draw (5.865450598678444,2.512891876085364)-- (5.652382758666241,2.512891876085364);
\draw (5.859090364648229,2.6814380778860616)-- (5.652382758666241,2.6782579608709542);
\draw (5.652382758666241,2.6782579608709542)-- (5.652382758666241,2.512891876085364);
\draw (6.,2.6)-- (5.865450598678444,2.512891876085364);
\draw (5.8,2.6)-- (5.652382758666241,2.512891876085364);
\draw (5.999118452553847,2.7728532194486912)-- (5.859090364648229,2.6814380778860616);
\draw (5.797172680635512,2.766733650602681)-- (5.652382758666241,2.6782579608709542);
\draw (5.706444747923069,2.4174883656321393) node[anchor=north west] {Figure 1.};
\begin{scriptsize}
\draw [fill=black] (5.999118452553847,2.7728532194486912) circle (1.5pt);
\draw[color=black] (6.013816683464034,2.788183109535685) node {$5$};
\draw [fill=black] (6.,2.6) circle (1.5pt);
\draw[color=black] (6.012996800479142,2.624195971614127) node {$1$};
\draw [fill=black] (5.8,2.6) circle (1.5pt);
\draw[color=black] (5.808208492406509,2.5796743334026218) node {$4$};
\draw [fill=black] (5.797172680635512,2.766733650602681) circle (1.5pt);
\draw[color=black] (5.785947673300757,2.788183109535685) node {$8$};
\draw [fill=black] (5.859090364648229,2.6814380778860616) circle (1.5pt);
\draw[color=black] (5.846369896587799,2.702319950127782) node {$6$};
\draw [fill=black] (5.865450598678444,2.512891876085364) circle (1.5pt);
\draw[color=black] (5.859090364648229,2.4874509399645044) node {$2$};
\draw [fill=black] (5.652382758666241,2.512891876085364) circle (1.5pt);
\draw[color=black] (5.6460225246360265,2.4810907059342893) node {$3$};
\draw [fill=black] (5.652382758666241,2.6782579608709542) circle (1.5pt);
\draw[color=black] (5.650581588515167,2.7036988969918148) node {$7$};
\end{scriptsize}
\end{tikzpicture}
\end{center}
Moreover, for a fixed integer $k\geq 2$ and for each integer $d$ satisfying
$1\leq d \leq 8\times 7^{k-2}$, we use the notation $Q_d^{(k)}$ to denote a cubic graph of
order~$8$. The vertex set of $Q_d^{(k)}$ is defined as
\[
V(Q_d^{(k)})=\{(x,1)_d^{(k)},\ldots,(x,8)_d^{(k)}\}.
\]
The edge set $E(Q_d^{(k)})$ is given by
\begin{align*}
	E(Q_d^{(k)})&=\{(x,i)_d^{(k)}(x,j)_d^{(k)} \mid i,j\in W_1,\ i<j,\ j-i\in\{1,3\}\}\\
	&\quad \cup \{(x,i)_d^{(k)}(x,j)_d^{(k)} \mid i,j\in W_2,\ i<j,\ j-i\in\{1,3\}\}\\
	&\quad \cup \{(x,i)_d^{(k)}(x,j)_d^{(k)} \mid i\in W_1,\ j\in W_2,\ j-i=4\}.
\end{align*}
The graph $Q_d^{(2)}$ is illustrated in Figure~2.

\begin{center}
\begin{tikzpicture}[line cap=round,line join=round,>=triangle 45,x=15.0cm,y=10.0cm]
\clip(5.531538312092156,2.3252649721940215) rectangle (6.113499725856829,2.865884864762297);
\draw (5.999118452553847,2.7728532194486912)-- (6.,2.6);
\draw (6.,2.6)-- (5.8,2.6);
\draw (5.8,2.6)-- (5.797172680635512,2.766733650602681);
\draw (5.797172680635512,2.766733650602681)-- (5.999118452553847,2.7728532194486912);
\draw (5.859090364648229,2.6814380778860616)-- (5.865450598678444,2.512891876085364);
\draw (5.865450598678444,2.512891876085364)-- (5.652382758666241,2.512891876085364);
\draw (5.859090364648229,2.6814380778860616)-- (5.652382758666241,2.6782579608709542);
\draw (5.652382758666241,2.6782579608709542)-- (5.652382758666241,2.512891876085364);
\draw (6.,2.6)-- (5.865450598678444,2.512891876085364);
\draw (5.8,2.6)-- (5.652382758666241,2.512891876085364);
\draw (5.999118452553847,2.7728532194486912)-- (5.859090364648229,2.6814380778860616);
\draw (5.797172680635512,2.766733650602681)-- (5.652382758666241,2.6782579608709542);
\draw (5.706444747923069,2.4174883656321393) node[anchor=north west] {Figure 2.};
\begin{scriptsize}
\draw [fill=black] (5.999118452553847,2.7728532194486912) circle (1.5pt);
\draw[color=black] (6.036816683464034,2.788183109535685) node {$(x, 5)_d^{(2)}$};
\draw [fill=black] (6.,2.6) circle (1.5pt);
\draw[color=black] (6.036996800479142,2.624195971614127) node {$(x, 1)_d^{(2)}$};
\draw [fill=black] (5.8,2.6) circle (1.5pt);
\draw[color=black] (5.758208492406509,2.6196743334026218) node {$(x, 4)_d^{(2)}$};
\draw [fill=black] (5.797172680635512,2.766733650602681) circle (1.5pt);
\draw[color=black] (5.785947673300757,2.798183109535685) node {$(x, 8)_d^{(2)}$};
\draw [fill=black] (5.859090364648229,2.6814380778860616) circle (1.5pt);
\draw[color=black] (5.846369896587799,2.712319950127782) node {$(x, 6)_d^{(2)}$};
\draw [fill=black] (5.865450598678444,2.512891876085364) circle (1.5pt);
\draw[color=black] (5.859090364648229,2.4874509399645044) node {$(x, 2)_d^{(2)}$};
\draw [fill=black] (5.652382758666241,2.512891876085364) circle (1.5pt);
\draw[color=black] (5.6460225246360265,2.4810907059342893) node {$(x, 3)_d^{(2)}$};
\draw [fill=black] (5.652382758666241,2.6782579608709542) circle (1.5pt);
\draw[color=black] (5.620581588515167,2.7036988969918148) node {$(x, 7)_d^{(2)}$};
\end{scriptsize}
\end{tikzpicture}
\end{center}
It is clear that the cubic graph $Q_d^{(k)}$ of order $8$, defined above, is isomorphic to the Cartesian product
$C_4\Box P_2$. Now consider the graph $CCC(n)$ with vertex set
\[
V(CCC(n)) = L_1 \cup L_2 \cup \cdots \cup L_n,
\]
where the sets $L_1, L_2, \ldots, L_n$ are referred to as the \emph{layers} of $CCC(n)$. The first layer is given by
$L_1 = V(C_4\Box P_2) = \{1,2,\ldots,8\}$. For a fixed integer $k\geq 2$ and each
$1\leq d \leq 8\times 7^{k-2}$, the cube $Q_d^{(k)}$, defined previously, is assumed to lie in the layer $L_k$.

Thus, the layer $L_1$ consists of the vertices $\{1,2,\ldots,8\}$, while the layers $L_2$ and $L_3$ consist of the cubes
$\{Q_1^{(2)}, Q_2^{(2)}, \ldots, Q_8^{(2)}\}$ and
$\{Q_1^{(3)}, Q_2^{(3)}, \ldots, Q_{56}^{(3)}\}$, respectively. In general, each layer
\[
L_k = \{Q_1^{(k)}, Q_2^{(k)}, \ldots, Q_{8\times 7^{k-2}}^{(k)}\},
\]
consists of $8\times 7^{k-2}$ cubes. Consequently, for $k\geq 2$, the layer $L_k$ contains
$64\times 7^{k-2}$ vertices. Therefore, the total number of vertices of $CCC(n)$ is
\[
|V(\text{CCC}(n))| = 8 + 64 \sum_{k=2}^{n} 7^{k-2}.
\]

We now describe the adjacency relation in the graph $CCC(n)$. Let $d$ be an arbitrary vertex in the layer $L_1$,
with $1\leq d\leq 8$. The vertex $d$ is adjacent to the vertex $(x,1)_d^{(2)}$ of the cube $Q_d^{(2)}$ in the layer $L_2$.
This construction yields the graph $CCC(2)$, as illustrated in Figure~3.

At the next level, each vertex of degree $3$ in every cube of the layer $L_2$ is adjacent, by an edge, to the vertex
$(x,1)_d^{(3)}$ of exactly one cube $Q_d^{(3)}$ in the layer $L_3$, resulting in the graph $CCC(3)$. More generally,
for each $k\geq 2$, every vertex of degree $3$ in each cube of the layer $L_k$ is adjacent to the vertex
$(x,1)_d^{(k+1)}$ of exactly one cube $Q_d^{(k+1)}$ in the layer $L_{k+1}$. This recursive procedure produces the graph
$CCC(k+1)$, and by continuing in this manner, the crystal cubic carbon graph $CCC(n)$ is obtained.

The layer
\[
L_n = \{Q_1^{(n)}, Q_2^{(n)}, \ldots, Q_{8\times 7^{n-2}}^{(n)}\},
\]
is called the \emph{outermost layer} of $CCC(n)$. Each vertex $(x,1)_d^{(n)}$ of degree $4$,
$1\leq d \leq 8\times 7^{n-2}$, in a cube $Q_d^{(n)}$ belonging to the outermost layer $L_n$ is referred to as the
\emph{head vertex} of $Q_d^{(n)}$.

\begin{center}
\begin{tikzpicture}[line cap=round,line join=round,>=triangle 45,x=10.0cm,y=10.0cm]
\clip(5.9081308513605615,1.2932437619909456) rectangle (7.099822839633116,2.101402696566582);
\draw (6.4,1.8)-- (6.5,1.8);
\draw (6.5,1.8)-- (6.5,1.7);
\draw (6.5,1.7)-- (6.4,1.7);
\draw (6.4,1.7)-- (6.4,1.8);
\draw (6.35,1.75)-- (6.45,1.75);
\draw (6.45,1.65)-- (6.35,1.65);
\draw (6.45,1.75)-- (6.45,1.65);
\draw (6.35,1.75)-- (6.35,1.65);
\draw (6.4,1.8)-- (6.35,1.75);
\draw (6.5,1.8)-- (6.45,1.75);
\draw (6.5,1.7)-- (6.45,1.65);
\draw (6.4,1.7)-- (6.35,1.65);
\draw (6.55,2.05)-- (6.65,2.05);
\draw (6.65,2.05)-- (6.65,1.95);
\draw (6.65,1.95)-- (6.55,1.95);
\draw (6.55,1.95)-- (6.55,2.05);
\draw (6.5,2.)-- (6.6,2.);
\draw (6.6,2.)-- (6.6,1.9);
\draw (6.6,1.9)-- (6.5,1.9);
\draw (6.5,1.9)-- (6.5,2.);
\draw (6.55,2.05)-- (6.5,2.);
\draw (6.65,2.05)-- (6.6,2.);
\draw (6.65,1.95)-- (6.6,1.9);
\draw (6.55,1.95)-- (6.5,1.9);
\draw (6.8,1.95)-- (6.95,1.95);
\draw (6.95,1.95)-- (6.95,1.85);
\draw (6.95,1.85)-- (6.8,1.85);
\draw (6.8,1.85)-- (6.8,1.95);
\draw (6.75,1.9)-- (6.9,1.9);
\draw (6.9,1.9)-- (6.9,1.8);
\draw (6.75,1.8)-- (6.75,1.9);
\draw (6.8,1.95)-- (6.75,1.9);
\draw (6.95,1.95)-- (6.9,1.9);
\draw (6.95,1.85)-- (6.9,1.8);
\draw (6.8,1.85)-- (6.75,1.8);
\draw (6.9,1.8)-- (6.75,1.8);
\draw (6.65,1.7)-- (6.75,1.7);
\draw (6.75,1.7)-- (6.75,1.6);
\draw (6.75,1.6)-- (6.65,1.6);
\draw (6.65,1.6)-- (6.65,1.7);
\draw (6.6,1.65)-- (6.7,1.65);
\draw (6.7,1.65)-- (6.7,1.55);
\draw (6.7,1.55)-- (6.6,1.55);
\draw (6.6,1.55)-- (6.6,1.65);
\draw (6.65,1.7)-- (6.6,1.65);
\draw (6.75,1.7)-- (6.7,1.65);
\draw (6.75,1.6)-- (6.7,1.55);
\draw (6.65,1.6)-- (6.6,1.55);
\draw (6.4,1.55)-- (6.5,1.55);
\draw (6.5,1.55)-- (6.5,1.45);
\draw (6.5,1.45)-- (6.4,1.45);
\draw (6.4,1.45)-- (6.4,1.55);
\draw (6.35,1.5)-- (6.45,1.5);
\draw (6.45,1.5)-- (6.45,1.4);
\draw (6.45,1.4)-- (6.35,1.4);
\draw (6.35,1.4)-- (6.35,1.5);
\draw (6.4,1.55)-- (6.35,1.5);
\draw (6.5,1.55)-- (6.45,1.5);
\draw (6.5,1.45)-- (6.45,1.4);
\draw (6.4,1.45)-- (6.35,1.4);
\draw (6.2,1.55)-- (6.3,1.55);
\draw (6.3,1.55)-- (6.3,1.45);
\draw (6.3,1.45)-- (6.2,1.45);
\draw (6.2,1.45)-- (6.2,1.55);
\draw (6.15,1.5)-- (6.25,1.5);
\draw (6.25,1.5)-- (6.25,1.4);
\draw (6.25,1.4)-- (6.15,1.4);
\draw (6.15,1.4)-- (6.15,1.5);
\draw (6.15,1.5)-- (6.2,1.55);
\draw (6.25,1.5)-- (6.3,1.55);
\draw (6.25,1.4)-- (6.3,1.45);
\draw (6.2,1.45)-- (6.15,1.4);
\draw (6.,1.7)-- (6.1,1.7);
\draw (6.1,1.7)-- (6.1,1.6);
\draw (6.1,1.6)-- (6.,1.6);
\draw (6.,1.6)-- (6.,1.7);
\draw (6.05,1.75)-- (6.15,1.75);
\draw (6.15,1.75)-- (6.15,1.65);
\draw (6.15,1.65)-- (6.05,1.65);
\draw (6.05,1.65)-- (6.05,1.75);
\draw (6.05,1.75)-- (6.,1.7);
\draw (6.15,1.75)-- (6.1,1.7);
\draw (6.15,1.65)-- (6.1,1.6);
\draw (6.05,1.65)-- (6.,1.6);
\draw (6.3,2.05)-- (6.4,2.05);
\draw (6.4,2.05)-- (6.4,1.95);
\draw (6.4,1.95)-- (6.3,1.95);
\draw (6.3,1.95)-- (6.3,2.05);
\draw (6.25,2.)-- (6.35,2.);
\draw (6.35,2.)-- (6.35,1.9);
\draw (6.35,1.9)-- (6.25,1.9);
\draw (6.25,1.9)-- (6.25,2.);
\draw (6.3,2.05)-- (6.25,2.);
\draw (6.4,2.05)-- (6.35,2.);
\draw (6.4,1.95)-- (6.35,1.9);
\draw (6.3,1.95)-- (6.25,1.9);
\draw (6.1,1.95)-- (6.2,1.95);
\draw (6.2,1.95)-- (6.2,1.85);
\draw (6.2,1.85)-- (6.1,1.85);
\draw (6.1,1.85)-- (6.1,1.95);
\draw (6.05,1.9)-- (6.15,1.9);
\draw (6.15,1.9)-- (6.15,1.8);
\draw (6.15,1.8)-- (6.05,1.8);
\draw (6.05,1.8)-- (6.05,1.9);
\draw (6.05,1.9)-- (6.1,1.95);
\draw (6.2,1.95)-- (6.15,1.9);
\draw (6.2,1.85)-- (6.15,1.8);
\draw (6.1,1.85)-- (6.05,1.8);
\draw (6.5,1.7)-- (6.75,1.8);
\draw (6.5,1.8)-- (6.5,1.9);
\draw (6.4,1.8)-- (6.35,1.9);
\draw (6.35,1.75)-- (6.2,1.85);
\draw (6.45,1.65)-- (6.5,1.55);
\draw (6.35,1.65)-- (6.3,1.55);
\draw (6.4,1.7)-- (6.15,1.65);
\draw (6.45,1.75)-- (6.6,1.65);
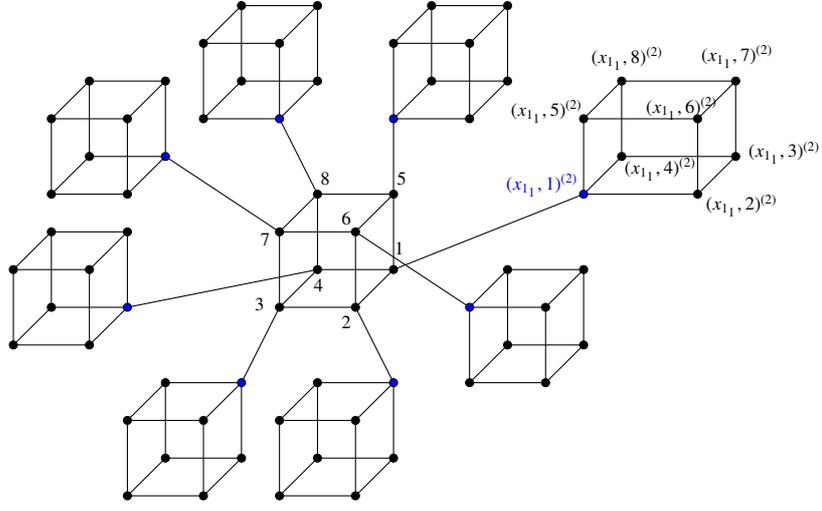
\draw (6.219059123178684,1.3573805022973364) node[anchor=north west] {Figure 3. Crystal cubic carbon $CCC(2)$};
\begin{scriptsize}
\draw [fill=black] (6.4,1.8) circle (1.5pt);
\draw[color=black] (6.4103765169160445,1.8213626891053447) node {$8$};
\draw [fill=black] (6.5,1.8) circle (1.5pt);
\draw[color=black] (6.520825650027141,1.8213626891053447) node {$5$};
\draw [fill=black] (6.5,1.7) circle (1.5pt);
\draw[color=black] (6.517781736902563,1.7285233388056953) node {$1$};
\draw [fill=black] (6.4,1.7) circle (1.5pt);
\draw[color=black] (6.401244777542308,1.6789524681861725) node {$4$};
\draw [fill=black] (6.35,1.75) circle (1.5pt);
\draw[color=black] (6.331234775676998,1.7398307306777459) node {$7$};
\draw [fill=black] (6.45,1.75) circle (1.5pt);
\draw[color=black] (6.437771735037253,1.7680942094252179) node {$6$};
\draw [fill=black] (6.45,1.65) circle (1.5pt);
\draw[color=black] (6.437771735037253,1.630249858192914) node {$2$};
\draw [fill=black] (6.35,1.65) circle (1.5pt);
\draw[color=black] (6.323624992865552,1.653947467253518) node {$3$};
\draw [fill=black] (6.55,2.05) circle (1.5pt);
\draw [fill=black] (6.65,2.05) circle (1.5pt);
\draw [fill=black] (6.65,1.95) circle (1.5pt);
\draw [fill=black] (6.55,1.95) circle (1.5pt);
\draw [fill=black] (6.5,2.) circle (1.5pt);
\draw [fill=black] (6.6,2.) circle (1.5pt);
\draw [fill=black] (6.6,1.9) circle (1.5pt);
\draw [fill=qqqqff] (6.5,1.9) circle (1.5pt);
\draw [fill=black] (6.8,1.95) circle (1.5pt);
\draw[color=black] (6.767607179673562,1.9796461715834353) node {$(x, 8)_1^{(2)}$};
\draw [fill=black] (6.95,1.95) circle (1.5pt);
\draw[color=black] (6.952193053091049,1.982690084708014) node {$(x, 7)_1^{(2)}$};
\draw [fill=black] (6.95,1.85) circle (1.5pt);
\draw[color=black] (7.01412618713943,1.85484573347571) node {$(x, 3)_1^{(2)}$};
\draw [fill=black] (6.8,1.85) circle (1.5pt);
\draw[color=black] (6.851080222178618,1.872977472849446) node {$(x, 4)_1^{(2)}$};
\draw [fill=black] (6.75,1.9) circle (1.5pt);
\draw[color=black] (6.701509351559095,1.9196361697181259) node {$(x, 5)_1^{(2)}$};
\draw [fill=black] (6.9,1.9) circle (1.5pt);
\draw[color=black] (6.869139138101161,1.924202039404994) node {$(x, 6)_1^{(2)}$};
\draw [fill=black] (6.9,1.8) circle (1.5pt);
\draw[color=black] (6.958056312784659,1.7848357316104007) node {$(x, 2)_1^{(2)}$};
\draw [fill=qqqqff] (6.75,1.8) circle (1.5pt);
\draw[color=qqqqff] (6.745421525309938,1.7722309497316087) node {$(x, 1)_1^{(2)}$};
\draw [fill=black] (6.65,1.7) circle (1.5pt);
\draw [fill=black] (6.75,1.7) circle (1.5pt);
\draw [fill=black] (6.75,1.6) circle (1.5pt);
\draw [fill=black] (6.65,1.6) circle (1.5pt);
\draw [fill=qqqqff] (6.6,1.65) circle (1.5pt);
\draw [fill=black] (6.7,1.65) circle (1.5pt);
\draw [fill=black] (6.7,1.55) circle (1.5pt);
\draw [fill=black] (6.6,1.55) circle (1.5pt);
\draw [fill=black] (6.4,1.55) circle (1.5pt);
\draw [fill=qqqqff] (6.5,1.55) circle (1.5pt);
\draw [fill=black] (6.5,1.45) circle (1.5pt);
\draw [fill=black] (6.4,1.45) circle (1.5pt);
\draw [fill=black] (6.35,1.5) circle (1.5pt);
\draw [fill=black] (6.45,1.5) circle (1.5pt);
\draw [fill=black] (6.45,1.4) circle (1.5pt);
\draw [fill=black] (6.35,1.4) circle (1.5pt);
\draw [fill=black] (6.2,1.55) circle (1.5pt);
\draw [fill=qqqqff] (6.3,1.55) circle (1.5pt);
\draw [fill=black] (6.3,1.45) circle (1.5pt);
\draw [fill=black] (6.2,1.45) circle (1.5pt);
\draw [fill=black] (6.15,1.5) circle (1.5pt);
\draw [fill=black] (6.25,1.5) circle (1.5pt);
\draw [fill=black] (6.25,1.4) circle (1.5pt);
\draw [fill=black] (6.15,1.4) circle (1.5pt);
\draw [fill=black] (6.,1.7) circle (1.5pt);
\draw [fill=black] (6.1,1.7) circle (1.5pt);
\draw [fill=black] (6.1,1.6) circle (1.5pt);
\draw [fill=black] (6.,1.6) circle (1.5pt);
\draw [fill=black] (6.05,1.75) circle (1.5pt);
\draw [fill=black] (6.15,1.75) circle (1.5pt);
\draw [fill=qqqqff] (6.15,1.65) circle (1.5pt);
\draw [fill=black] (6.05,1.65) circle (1.5pt);
\draw [fill=black] (6.3,2.05) circle (1.5pt);
\draw [fill=black] (6.4,2.05) circle (1.5pt);
\draw [fill=black] (6.4,1.95) circle (1.5pt);
\draw [fill=black] (6.3,1.95) circle (1.5pt);
\draw [fill=black] (6.25,2.) circle (1.5pt);
\draw [fill=black] (6.35,2.) circle (1.5pt);
\draw [fill=qqqqff] (6.35,1.9) circle (1.5pt);
\draw [fill=black] (6.25,1.9) circle (1.5pt);
\draw [fill=black] (6.1,1.95) circle (1.5pt);
\draw [fill=black] (6.2,1.95) circle (1.5pt);
\draw [fill=qqqqff] (6.2,1.85) circle (1.5pt);
\draw [fill=black] (6.1,1.85) circle (1.5pt);
\draw [fill=black] (6.05,1.9) circle (1.5pt);
\draw [fill=black] (6.15,1.9) circle (1.5pt);
\draw [fill=black] (6.15,1.8) circle (1.5pt);
\draw [fill=black] (6.05,1.8) circle (1.5pt);
\end{scriptsize}
\end{tikzpicture}
\end{center}
\section{Results for the crystal cubic carbon $CCC(n)$}
We start this section with the following theorem:  
\begin{theorem} {\rm\cite{Zhang-1}} \label{c.1}
	Let $CCC(n)$ be the crystal cubic carbon graph. If $n\geq 2$ is an integer, then the minimum
	cardinality of a resolving set of $CCC(n)$ is $16\times 7^{\,n-2}$.
\end{theorem}
\begin{observation}\label{c.2}
	Let $S$ be an arbitrary minimal resolving set of $CCC(n)$. According to the proof of
	Theorem~\ref{c.1}, the set $S$ contains exactly two vertices from each cube in the outermost
	layer
	\[
	L_n=\{Q_1^{(n)}, Q_2^{(n)}, \ldots, Q_{8\times 7^{n-2}}^{(n)}\},
	\]
	of $CCC(n)$. Thus, if
	$S=\{S_1,S_2,\ldots,S_{8\times 7^{n-2}}\}$ is a minimal resolving set such that
	$S_d\subset Q_d^{(n)}$, then each $S_d$ consists of two vertices from the cube $Q_d^{(n)}$.
\end{observation}
	We now determine which vertices of each cube $Q_d^{(n)}$ in the outermost layer $L_n$ may
	belong to $S_d$. Without loss of generality, consider the cube $Q_1^{(n)}$ and let
	$x=(x,i)_d^{(n)}$ be an arbitrary vertex from another cube $Q_d^{(n)}$ in $L_n$,
	$d\neq 1$, such that $x\in S_d$. Suppose that the distance between the head vertex
	$(x,1)_1^{(n)}\in Q_1^{(n)}$ and $x$ is a positive integer $c$, that is,
	$r((x,1)_1^{(n)}\mid x)=c$. Then any two vertices $(x,i)_1^{(n)}$ and $(x,j)_1^{(n)}$ of
	$Q_1^{(n)}$ may belong to $S_1$, provided that the set
	\[
	\{(x,i)_1^{(n)}, (x,j)_1^{(n)}, (x,i)_d^{(n)}\},
	\]
	resolves all vertices of $Q_1^{(n)}$. Consequently, $S_1$ may consist of any of the following
	pairs of vertices of $Q_1^{(n)}$:
	\[
	\begin{aligned}
	&(x,2)_1^{(n)}, (x,3)_1^{(n)};\ (x,2)_1^{(n)}, (x,4)_1^{(n)};\ (x,2)_1^{(n)}, (x,5)_1^{(n)};\\
	&(x,2)_1^{(n)}, (x,6)_1^{(n)};\ (x,2)_1^{(n)}, (x,8)_1^{(n)};\\
	&(x,3)_1^{(n)}, (x,4)_1^{(n)};\ (x,3)_1^{(n)}, (x,6)_1^{(n)};\ (x,3)_1^{(n)}, (x,8)_1^{(n)};\\
	&(x,4)_1^{(n)}, (x,5)_1^{(n)};\ (x,4)_1^{(n)}, (x,8)_1^{(n)};\\
	&(x,5)_1^{(n)}, (x,6)_1^{(n)};\ (x,5)_1^{(n)}, (x,8)_1^{(n)};\\
	&(x,6)_1^{(n)}, (x,8)_1^{(n)}.
	\end{aligned}
	\]
\begin{proposition}\label{c.3}
	Let $S$ be a minimal resolving set of $CCC(n)$. Then the head vertex $(x,1)_d^{(n)}$ of any
	cube in the outermost layer $L_n$ does not belong to $S$.
\end{proposition}

\begin{proof}
	Let $S=\{S_1,S_2,\ldots,S_{8\times 7^{n-2}}\}$ be a minimal resolving set of $CCC(n)$ with
	$S_d\subset Q_d^{(n)}$. As noted earlier, each $S_d$ contains exactly two vertices from
	$Q_d^{(n)}$. Suppose, to the contrary, that for the cube $Q_1^{(n)}$ the set $S_1$ contains
	the head vertex $(x,1)_1^{(n)}$.
	
	Let $x=(x,i)_d^{(n)}\in S_d$ be an arbitrary vertex from another cube
	$Q_d^{(n)}\subset L_n$, $d\neq 1$, and assume that
	$r((x,1)_1^{(n)}\mid x)=c$ for some positive integer $c$. Then the set
	\[
	\{(x,1)_1^{(n)}, (x,j)_1^{(n)}, (x,i)_d^{(n)}\},
	\]
	fails to resolve all vertices of $Q_1^{(n)}$. Hence,
	$S=\{S_1,S_2,\ldots,S_{8\times 7^{n-2}}\}$ cannot be a resolving set of $CCC(n)$, which
	contradicts the minimality of $S$. Therefore, no head vertex belongs to a minimal resolving
	set.
\end{proof}
\begin{proposition}\label{c.4}
	Let $S$ be a minimal resolving set of $CCC(n)$. Then the vertex $(x,7)_d^{(n)}$ of any cube in
	the outermost layer $L_n$ does not belong to $S$.
\end{proposition}

\begin{proof}
	The argument follows analogously to the proof of Proposition~\ref{c.3}.
\end{proof}
\begin{proposition}\label{c.5}
	Let $W$ be a resolving set of $CCC(n)$ consisting of exactly two vertices adjacent to the head
	vertex in the outermost layer $L_n$ of each cube in $L_n$. Then $W$ is a minimal resolving set
	of $CCC(n)$.
\end{proposition}

\begin{proof}
	Recall that $V(CCC(n))=L_1\cup L_2\cup \cdots \cup L_n$, where $L_1,L_2,\ldots,L_n$ are the
	layers of $CCC(n)$. Define
	\[
	Z_1=\{(x,2)_1^{(n)}, (x,2)_2^{(n)}, \ldots, (x,2)_{8\times 7^{n-2}}^{(n)}\},
	\]
	and
	\[
	Z_2=\{(x,4)_1^{(n)}, (x,4)_2^{(n)}, \ldots, (x,4)_{8\times 7^{n-2}}^{(n)}\}.
	\]
	Each of the sets $Z_1$ and $Z_2$ contains exactly one vertex adjacent to the head vertex from
	each cube in the outermost layer $L_n$. By Observation~\ref{c.2}, the set
	$W=Z_1\cup Z_2$, which consists of exactly two such adjacent vertices from each cube of $L_n$,
	forms a minimal resolving set of $CCC(n)$.
\end{proof}

\begin{conclusion}\label{c.6}
	 Let  
	$S=\{S_1, S_2, \ldots, S_{8\times 7^{n-2}}\}$ be an arbitrary minimal resolving set of $CCC(n)$ such that  
	$S_d \subset Q_d^{(n)}$. Then each $S_d$ contains any pair of elements described in Observation~\ref{c.2} from each cube $Q_d^{(n)}$ in the outermost layer $L_n$ of $CCC(n)$.
\end{conclusion}


\begin{proposition}\label{c.7}
 If  
	$S=\{S_1, S_2, \ldots, S_{8\times 7^{n-2}}\}$ is an arbitrary minimal resolving set of $CCC(n)$ such that  
	$S_d \subset Q_d^{(n)}$, then $S$ cannot be a doubly resolving set of $CCC(n)$.
\end{proposition}

\begin{proof}
	Let
	\[
	V(CCC(n)) = L_1 \cup L_2 \cup \cdots \cup L_n,
	\]
	where the sets $L_1, L_2, \ldots, L_n$ are called the layers of $CCC(n)$, with
	\[
	L_1 = V(C_4 \square P_2) = \{1,2,\ldots,8\},
	\]
	and for $k \geq 2$,
	\[
	L_k = \{Q_1^{(k)}, Q_2^{(k)}, \ldots, Q_{8\times 7^{k-2}}^{(k)}\},
	\]
	where each $Q_d^{(k)}$ denotes a cube in the layer $L_k$.
	
	By Conclusion~\ref{c.6}, if  
	$S=\{S_1, S_2, \ldots, S_{8\times 7^{n-2}}\}$ is a minimal resolving set of $CCC(n)$ such that  
	$S_d \subset Q_d^{(n)}$, then each $S_d$ contains any pair of elements defined in Observation~\ref{c.2} from each cube $Q_d^{(n)}$ in the outermost layer
	\[
	L_n=\{Q_1^{(n)}, Q_2^{(n)}, \ldots, Q_{8\times 7^{n-2}}^{(n)}\}.
	\]
	
	Without loss of generality, consider the cube $Q_1^{(n)}$ in the outermost layer $L_n$.
	Let $(x,i)_1^{(n)}$ and $(x,j)_1^{(n)}$ be two arbitrary vertices of $Q_1^{(n)}$ such that
	$S_1$ contains both $(x,i)_1^{(n)}$ and $(x,j)_1^{(n)}$.
	Let $x=(x,i)_d^{(n)}$ be an arbitrary vertex belonging to another cube $Q_d^{(n)}$ in the outermost layer $L_n$, with $d \neq 1$, such that $x \in S_d$.
	Assume that the distance between the head vertex $(x,1)_1^{(n)} \in Q_1^{(n)}$ and $x$ is a positive integer $c$, that is,
	\[
	r\big((x,1)_1^{(n)} \mid x\big)=c.
	\]
	
	Let
	\[
	Z=\{(x,i)_1^{(n)}, (x,j)_1^{(n)}, x\} \subseteq S.
	\]
	Then the set $Z$ cannot doubly resolve all vertices of the cube $Q_1^{(n)}$.
	Hence, if $S=\{S_1, S_2, \ldots, S_{8\times 7^{n-2}}\}$ is an arbitrary minimal resolving set of $CCC(n)$ such that
	$S_d \subset Q_d^{(n)}$, then $S$ cannot be a doubly resolving set of $CCC(n)$.
\end{proof}


\begin{conclusion}\label{c.8}
 If $S$ is a doubly resolving set of $CCC(n)$, then at least three vertices from each cube in the outermost layer $L_n$ of $CCC(n)$ must belong to $S$.
\end{conclusion}


\begin{theorem}\label{c.9}
	If $n \geq 2$ is an integer, then the minimum cardinality of a doubly resolving set of $CCC(n)$ is
	\[
	3 \times 8 \times 7^{n-2}.
	\]
\end{theorem}

\begin{proof}
	Let
	\[
	V(CCC(n)) = L_1 \cup L_2 \cup \cdots \cup L_n,
	\]
	be the vertex set of $CCC(n)$, where the layers $L_1, L_2, \ldots, L_n$ are defined as above.
	By Conclusion~\ref{c.8}, if $S$ is a doubly resolving set of $CCC(n)$, then at least three vertices from each cube in the outermost layer $L_n$ must belong to $S$.
	Therefore,
	\[
	|S| \geq 3 \times 8 \times 7^{n-2}.
	\]
	
	Now, let
	\[
	S=\{S_1, S_2, \ldots, S_{8\times 7^{n-2}}\},
	\]
	be a set such that each $S_d \subset Q_d^{(n)}$ consists of exactly three vertices adjacent to the head vertex of the cube $Q_d^{(n)}$ in the outermost layer $L_n$.
	We show that $S$ is a minimal doubly resolving set of $CCC(n)$.
	
	Consider the cube $Q_1^{(n)}$ in the outermost layer
	\[
	L_n=\{Q_1^{(n)}, Q_2^{(n)}, \ldots, Q_{8\times 7^{n-2}}^{(n)}\}.
	\]
	Let the three vertices
	\[
	(x,2)_1^{(n)},\ (x,4)_1^{(n)},\ (x,5)_1^{(n)},
	\]
	belong to $S_1$.
	Let
	\[
	Z=\{(x,2)_1^{(n)}, (x,4)_1^{(n)}, (x,5)_1^{(n)}, x\},
	\]
	where $x=(x,i)_d^{(n)} \in S_d$ is a vertex from another cube in $L_n$ with $d \neq 1$, and
	\[
	r\big((x,1)_1^{(n)} \mid x\big)=c
	\]
	for some positive integer $c$.
	Then we obtain:
	\[
	\begin{aligned}
	r((x,1)_1^{(n)} \mid Z) &= (1,1,1,c),\\
	r((x,2)_1^{(n)} \mid Z) &= (0,2,2,c+1),\\
	r((x,3)_1^{(n)} \mid Z) &= (1,1,3,c+2),\\
	r((x,4)_1^{(n)} \mid Z) &= (2,0,2,c+1),\\
	r((x,5)_1^{(n)} \mid Z) &= (2,2,0,c+1),\\
	r((x,6)_1^{(n)} \mid Z) &= (1,3,1,c+2),\\
	r((x,7)_1^{(n)} \mid Z) &= (2,2,2,c+3),\\
	r((x,8)_1^{(n)} \mid Z) &= (3,1,1,c+2).
	\end{aligned}
	\]
	
	Hence, every vertex of $Q_1^{(n)}$ is doubly resolved by the set $Z$.
	Therefore, if each $S_d$ contains three adjacent vertices with respect to the head vertex of each cube in the outermost layer $L_n$, then
	\[
	S=\{S_1, S_2, \ldots, S_{8\times 7^{n-2}}\},
	\]
	is a minimal doubly resolving set of $CCC(n)$.
	Thus, the minimum cardinality of a doubly resolving set of $CCC(n)$ is
	\[
	3 \times 8 \times 7^{n-2}.
	\]
\end{proof}

\begin{proposition}\label{c.10}
	 If  
	$S=\{S_1, S_2, \ldots, S_{8\times 7^{n-2}}\}$ is an arbitrary minimal resolving set of $CCC(n)$ such that  
	$S_d \subset Q_d^{(n)}$, then $S$ cannot be a strong resolving set of $CCC(n)$.
\end{proposition}

\begin{proof}
	Let
	\[
	V(CCC(n)) = L_1 \cup L_2 \cup \cdots \cup L_n,
	\]
	where the sets $L_1, L_2, \ldots, L_n$ are called the layers of $CCC(n)$, with
	\[
	L_1 = V(C_4 \square P_2) = \{1,2,\ldots,8\},
	\]
	and for $k \geq 2$,
	\[
	L_k = \{Q_1^{(k)}, Q_2^{(k)}, \ldots, Q_{8\times 7^{k-2}}^{(k)}\},
	\]
	where each $Q_d^{(k)}$ denotes a cube in the layer $L_k$.
	
	By Conclusion~\ref{c.6}, if
	$S=\{S_1, S_2, \ldots, S_{8\times 7^{n-2}}\}$ is an arbitrary minimal resolving set of $CCC(n)$ such that
	$S_d \subset Q_d^{(n)}$, then each $S_d$ contains any pair of vertices described in Observation~\ref{c.2} from each cube $Q_d^{(n)}$ in the outermost layer
	\[
	L_n = \{Q_1^{(n)}, Q_2^{(n)}, \ldots, Q_{8\times 7^{n-2}}^{(n)}\}.
	\]
	
	Without loss of generality, consider an arbitrary cube $Q_d^{(n)}$ in the outermost layer $L_n$.
	For each such cube, there exist exactly two vertices, say $(x,i)_d^{(n)}$ and $(x,j)_d^{(n)}$, distinct from
	$(x,1)_d^{(n)}$ and $(x,7)_d^{(n)}$, which are mutually maximally distant in $Q_d^{(n)}$.
	Hence, at least one of these vertices must belong to any minimal strong resolving set of $CCC(n)$.
	
	Consequently, the minimum cardinality of a strong resolving set of $CCC(n)$ is at least
	\[
	3 \times 8 \times 7^{n-2}.
	\]
	Therefore, if $S=\{S_1, S_2, \ldots, S_{8\times 7^{n-2}}\}$ is an arbitrary minimal resolving set of $CCC(n)$ such that
	$S_d \subset Q_d^{(n)}$, then $S$ cannot be a strong resolving set of $CCC(n)$.
\end{proof}


\begin{proposition}\label{c.11}
	 If $U$ is an arbitrary strong resolving set of $CCC(n)$, then all but one of the vertices
	$(x,7)_d^{(n)}$ from the cubes $Q_d^{(n)}$ in the outermost layer $L_n$ of $CCC(n)$ must belong to $U$.
\end{proposition}

\begin{proof}
	By Proposition~\ref{c.4}, if $S$ is an arbitrary minimal resolving set of $CCC(n)$, then the vertices
	$(x,7)_d^{(n)}$ from each cube $Q_d^{(n)}$ in the outermost layer $L_n$ do not belong to $S$.
	The total number of such vertices is $8 \times 7^{n-2}$.
	Moreover, any pair of these vertices is mutually maximally distant in $CCC(n)$.
	
	Therefore, in each cube $Q_d^{(n)}$ of the outermost layer $L_n$, all but one of the vertices $(x,7)_d^{(n)}$
	must belong to any strong resolving set of $CCC(n)$.
	Hence, the number of vertices of this type contained in an arbitrary strong resolving set $U$ is
	\[
	8 \times 7^{n-2} - 1.
	\]
\end{proof}


\begin{theorem}\label{c.12}
	If $n \geq 2$ is an integer, then the minimum cardinality of a strong resolving set of $CCC(n)$ is
	\[
	32 \times 7^{n-2} - 1.
	\]
\end{theorem}

\begin{proof}
	Let
	\[
	V(CCC(n)) = L_1 \cup L_2 \cup \cdots \cup L_n,
	\]
	be the vertex set of $CCC(n)$, where the layers $L_1, L_2, \ldots, L_n$ are defined as above.
	Let
	\[
	S=\{S_1, S_2, \ldots, S_{8\times 7^{n-2}}\},
	\]
	be an arbitrary minimal resolving set of $CCC(n)$ such that $S_d \subset Q_d^{(n)}$.
	By Theorem~\ref{c.1}, each $S_d$ contains any pair of vertices described in Observation~\ref{c.2} from each cube $Q_d^{(n)}$ in the outermost layer
	\[
	L_n = \{Q_1^{(n)}, Q_2^{(n)}, \ldots, Q_{8\times 7^{n-2}}^{(n)}\}.
	\]
	
	It can be verified that any pair of vertices lying in the layers
	$L_1 \cup L_2 \cup \cdots \cup L_{n-1}$ is strongly resolved by some vertex of $S$.
	Thus, it remains to consider the vertices in the outermost layer $L_n$.
	
	By Proposition~\ref{c.10}, the minimum cardinality of a strong resolving set of $CCC(n)$ is at least
	\[
	3 \times 8 \times 7^{n-2}.
	\]
	On the other hand, by Proposition~\ref{c.4}, the vertices $(x,7)_d^{(n)}$ from each cube $Q_d^{(n)}$ in the outermost layer $L_n$ do not belong to $S$.
	There are $8 \times 7^{n-2}$ such vertices, and any two of them are mutually maximally distant in $CCC(n)$.
	Therefore, by Proposition~\ref{c.11}, all but one of these vertices must belong to any strong resolving set of $CCC(n)$.
	Hence, the minimum cardinality of a strong resolving set of $CCC(n)$ is
	\[
	32 \times 7^{n-2} - 1.
	\]
\end{proof}


\section{Conclusion}
In this work, we introduced an alternative construction of the crystal cubic carbon graph $CCC(n)$.
We determined the minimum cardinalities of doubly resolving sets and strong resolving sets for this family of graphs.
Further investigations may focus on computing various topological indices and exploring additional structural properties of $CCC(n)$.

\section*{Funding and Conflict of interest} 
The authors have not disclosed any funding and declare no conflict of interest.\\

\noindent \textbf{Acknowledgements}\\
We would like to express our gratitude to Prof. Jia-Bao Liu, who gave us help and insightful suggestions  provided a great improvement to our paper.  \\[2mm]


\bigskip
\end{document}